# A Multi-objective Sequential Quadratic Programming Algorithm Based on Low-order Smooth Penalty Function

Zanyang Kong

College of Mathematics and Statistics, Chongqing University,

Email:202206021062t@stu.cqu.edu.cn

**Abstract:** In this paper, we propose a Multi-Objective Sequential Quadratic Programming (MOSQP) algorithm for constrained multi-objective optimization problems, basd on a low-order smooth penalty function as the merit function for line search. The algorithm constructs single-objective optimization subproblems based on each objective function, solves quadratic programming (QP) subproblems to obtain descent directions for expanding the iterative point set within the feasible region, and filters non-dominated points after expansion. A new QP problem is then formulated using information from all objective functions to derive descent directions. The Armijo step size rule is employed for line search, combined with Powell's correction formula (1978) for $B^j$ iteration updates. If QP subproblems is infesible, the negative gradient of the merit function is adopted as the search direction. The algorithm is proven to converge to an approximate Pareto front for constrained multi-objective optimization. Finally, numerical experiments are performed for specific multi-objective optimization problems

# 0 Introduction

In 2015, Fliege and Vaz [1] using the idea of single target SQP method, developed the SQP method to solve the constraint multi-objective optimization problem, the algorithm is divided into two stages, by solving the secondary planning subproblem obtain search direction, it uses the L1 accurate penalty function as a value function line search step, to obtain approximate Pareto frontier with constraint multi-objective optimization problem. Fliege and Vaz discuss the problem of its convergence to a local Pareto optimal solution under appropriate differential assumptions. Extensive numerical experiments also confirm the superiority of this algorithm compared to the state-of-the-art multi-objective optimization solver NSGA-II and the classical scalarization method MOScalar. However, the value function in the paper uses the most classical L1 exact penalty function, and only the symmetric positive definite matrix is the Hesse matrix or the unit matrix is discussed. Later, we can consider correcting the symmetric positive definite matrix, and avoid a large calculation of Hesse matrix, and avoid the disadvantages such as using only the unit matrix. And the sub-problems in the algorithm are not necessarily feasible at every iterative step. In 2019, Gebken et al. [2] also discussed the solution of unconstrained and constrained multi-objective optimization problems by solving quadratic planning subproblems for multi-objective optimization problems. The iterative process in the paper does not use the penalty function, and only positive constraints are used in the subproblems. However, the SQP method in the paper requires a feasible initial approximation, which is very difficult to implement in the nonlinear constraint problem. In 2021, Ansary, Panda [3] [4] proposed a SQP algorithm different from [1] by proposing different secondary planning sub-problems and penalty functions, and had good computing performance compared with the original algorithm and classical weights and methods.

In this paper, we propose a new MOSQP algorithm based on the multi-objective optimization problems with inequality constraints, and prove that it converges to the Pareto critical point of the original multi-objective optimization problem under certain conditions. Finally, numerical experiments are conducted on multi-objective optimization problems with boundary constraints or nonlinear constraints, and the algorithm performance is compared by NSGAII and MOSQP(F).

# 1 Preliminaries

## 1.1 Basic definition

**Gradient:** let $f: R^n \to R$, and $f(x)$ is a continuous differentiable function, then the first derivative is defined as:

$$\nabla_x f(x) = (\frac{\partial f(x)}{\partial x_1}, \frac{\partial f(x)}{\partial x_2}, \cdots \frac{\partial f(x)}{\partial x_n})^T$$

**Hessian Matrix:** let $f: R^n \to R$, and $f(x)$ is a second-order continuous differentiable function, then say:

$$\nabla_{xx} f(x) = \begin{bmatrix} \frac{\partial^2 f(x)}{\partial x_1^2} & \frac{\partial^2 f(x)}{\partial x_1 \partial x_2} & \cdots & \frac{\partial^2 f(x)}{\partial x_1 \partial x_n} \\ \frac{\partial^2 f(x)}{\partial x_2 \partial x_1} & \frac{\partial^2 f(x)}{\partial x_2^2} & \cdots & \frac{\partial^2 f(x)}{\partial x_2 \partial x_n} \\ \vdots & \vdots & \ddots & \vdots \\ \frac{\partial^2 f(x)}{\partial x_n \partial x_1} & \frac{\partial^2 f(x)}{\partial x_n \partial x_2} & \cdots & \frac{\partial^2 f(x)}{\partial x_n^2} \end{bmatrix},$$

is the Hesse matrix of $f(x)$ at $x = (x_1, x_2, \cdots x_n)^T \in X \subset R^n$.

**Jacobi Matrix**: let $F : R^n \to R^m$, then say:

$$J_F(x_1, x_2, \cdots x_n) = \begin{bmatrix} \frac{\partial f_1(x)}{\partial x_1} & \frac{\partial f_1(x)}{\partial x_2} & \cdots & \frac{\partial f_1(x)}{\partial x_n} \\ \frac{\partial f_2(x)}{\partial x_1} & \frac{\partial f_2(x)}{\partial x_2} & \cdots & \frac{\partial f_2(x)}{\partial x_n} \\ \vdots & \vdots & \ddots & \vdots \\ \frac{\partial f_m(x)}{\partial x_1} & \frac{\partial f_m(x)}{\partial x_2} & \cdots & \frac{\partial f_m(x)}{\partial x_n} \end{bmatrix},$$

is the Jacobi matrix of $F(x) = (f_1(x), \cdots, f_m(x))^T, m \geq 2$ at $x = (x_1, x_2, \cdots x_n)^T \in X \subset R^n$

## 1.2 Problem

This paper focuses on the following multi-objective optimization problems of inequality constraints：

$$\text{(CMOP)} \quad \min \ F(x) = (f_1(x), \cdots, f_m(x))^T, m \geq 2$$
$$\text{s.t.} \quad g_i(x) \leq 0, i \in I = \{1, 2, \cdots a\}$$

where $x = (x_1, x_2, x_3, \cdots, x_n)^T \in X \subset R^n$, $f_i(x) : R^n \to R$, $g_i(x) : R^n \to R$ are second-order continuous differentiable functions, $X = \{x | g_i(x) \leq 0, i \in I\}$ is the feasible domain of the problem. Define the following symbols:

$$I_0(x) = \{i \in I | g_i(x) = 0\}, \quad I_+(x) = \{i \in I | g_i(x) > 0\}, \quad I_-(x) = \{i \in I | g_i(x) < 0\}$$

**Globally optimal solution:** If $x^* = (x_1, x_2, \cdots x_n)^T \in X \subset R^n$, for each $i = 1, \cdots, m, x \in X$ there are:

$$f_i(x^*) \leq f_i(x)$$

Then $x^*$ is the global optimal solution of problem (CMOP).

**Pareto Optimal Solution:** Let $x^* \in X$, If there is not $x \neq x^* \in X$, s.t. $f_i(x) \leq f_i(x^*), (i=1,2,\cdots m)$ and at least one $i_0, s.t. f_{i_0}(x) < f_{i_0}(x^*)$, it is called a Pareto optimal solution of the problem (CMOP).

**Pareto Front:** Let $X^*$ be the full optimal solution of the problem, then $F(X^*) = \{F(x) | x \in X^*\}$ is called as the Pareto Front of the multi-objective optimization problem.

**KKT condition:** If $x^* \in X \subset R^n$, $\lambda_i^* > 0, i \in I$, $w_i^* \geq 0$ and at least one $w_i^* > 0$, and the following conditions:

$$\sum_{i=1}^{m} w_i^* \nabla g_i(x^*) + \sum_{i=1}^{a} \lambda_i^* \nabla g_i(x^*) = 0$$

$$g_i(x^*) \leq 0, \lambda_i^* > 0, \forall i \in I$$

$$\lambda_i^* g_i(x^*) = 0, \forall i \in I$$

$x^* \in X \subset R^n$ is called the KKT point of the original multi-objective optimization problem (CMOP).

**Pareto critical point:**
If $x^*$ is the KKT point of the problem:

$$\min \sum_{i=1}^{m} w_i f_i(x)$$

$$s.t. g_i(x) \leq 0, i \in I$$

The KKT point of, where $w_i \geq 0$, and at least one $w_i > 0$. Then $x^*$ is called the Pareto critical point of the original multi-objective constraint optimization problem (CMOP).

## 1.3 MOSQP Algorithm

In 2015, Fliege and Vaz [1] using the idea of single target SQP method, developed the SQP method to solve the constraint multi-objective optimization problem, the algorithm is divided into two stages, by solving the secondary planning subproblem obtain search direction, it uses the L1 accurate penalty function as a value function line search step, to obtain approximate Pareto frontier with constraint multi-objective optimization problem. The algorithmic framework is provided as follows:

**Algorithm 2:**
**Step1** Initialize the iteration point set $X_0 = \{x_1, x_2, \cdots x_N\}, x_i \in R^n$;
**Step2** Spread out K times for each iteration point in $X_j = \{x_1^j, x_2^j, \cdots x_N^j\}, x_i^j \in R^n$.

Let $T = \varnothing$,

Consider each objective function to construct the corresponding single-objective optimization problem:

$$(\text{P}_i) \quad \min \ f_i(x), i \in M = \{1, \cdots, m\}$$
$$\text{s.t.} \quad g_i(x) \leq 0, i \in I = \{1, 2, \cdots a\}$$

Solving for the quadratic programming sub-problem

$$\mathbf{QP}(x^j, B^j) \quad \min \frac{1}{2}(d^j)^T B^j d^j + \nabla f_i(x^j)^T d^j$$
$$\text{s.t.} \ \nabla g_i(x^j)^T d + g_i(x^j) \leq 0 \quad i \in I$$

And a linear search by the Armijo step rule:

$P(x^j + \alpha^j d^j, \pi) \leq P(x^j, \pi) + \sigma \alpha^j \nabla P(x^j, \pi) d^j$, where $\sigma \in (0,1)$。

Let $T = T \cup \{x^j + \alpha^j d^j\}$

Let $X^{j+1} = \text{non-dominate}\{X^j \cup T\}$, $j = j+1$

If $j = K$, stop, let $X = X^K$.

**Step3** Iteratively optimization for each point set in $X = \{x_1^0, x_2^0, \cdots, x_{M^*}^0\}$.

Let $S = \varnothing$:

Obtaining the iteration direction by solving a quadratic programming subproblem with information about all the objective functions:

$$\mathbf{QP}(x^j, B^j) \quad \min \frac{1}{2}(d^j)^T B^j d^j + \sum_{i=1}^{m} \nabla f_i(x^j)^T d^j$$
$$\text{s.t.} \ \nabla g_i(x^j)^T d^j + g_i(x^j) \leq 0 \quad i \in I$$

And a linear search by the Armijo step rule:

$P(x^j + \alpha^j d^j, \pi) \leq P(x^j, \pi) + \sigma \alpha^j \nabla P(x^j, \pi) d^j$, where $\sigma \in (0,1)$。

Let $x^{j+1} = x^j + \alpha^j d^j$, $j = j+1$, if $\|d^j\| = 0$, Stop。

$S = S \cup x_l^j, l = \{1, 2, \cdots, M^*\}$, Iterate over the next iteration point.

The non-dominant point in S is retained as the final Pareto optimal solution.

In Spread stage, this paper proposes a feasible spread mode to ensure that the diffusion stage is only carried out in the feasible domain, so that the initial iteration point set of the Pareto optimization solution stage is all in the feasible domain, and we cancel the non-dominant solution at each diffusion, which is because some non-dominant solutions still retain the possibility of declining towards a certain function.

# 2 A smooth approximation of the low-order penalty functions

## 2.1 Smooth function

Consider the L1-type penalty function

$$P(x,\pi) = f(x) + \pi \sum_{i=1}^{l} g_i^+(x)$$

Introducing a function $h: R \to R$

$$h(t) = \begin{cases} 0 & t < 0 \\ t & t \geq 0 \end{cases}$$

Then

$$g_i^+(x) = h(g_i(x))$$

According to the function $h(t)$, then

$$h^k(t) = \begin{cases} 0 & t < 0 \\ t^k & t \geq 0 \end{cases}$$

When $k \in (0,1]$, it cannot be able Thus increasing the difficulty for considering a low-order penalty function as a value function. Consider the following smoothness approximation function:

$$h_\varepsilon^k(t) = \begin{cases} 0 & t < -\varepsilon \\ \varepsilon^{k(1-b)} b^{-k}(t+\varepsilon)^{kb} & -\varepsilon \leq t < 0 \\ (t + \varepsilon b^{-1})^k & t \geq 0 \end{cases}$$

Where $\varepsilon > 0$, $b > 0$.

## 2.2 Properties of the smooth functions

**Theorem 2.2.1** For any $0 < k < +\infty$, $\varepsilon > 0$, $b > 0$, the following properties to hold true:

(1) $h_\varepsilon^k(t)$ is continuous in $R$;

(2) $\forall t \in R, h_\varepsilon^m(t) \geq h^m(t)$;

(3) $\forall t \in R, \lim_{\varepsilon \to 0} h_\varepsilon^m(t) = h^m(t)$;

(4) If $0 < \varepsilon_1 < \varepsilon_2$, then $\forall t \in R, h_{\varepsilon_1}^k(t) \leq h_{\varepsilon_2}^k(t)$。

Proof：

(1) Just prove that $h_\varepsilon^k(t)$ is continuous in $t=0$, $t=-\varepsilon$

when $t=0$,
$$\lim_{t\to 0^-} h_\varepsilon^k(t) = \varepsilon^{k(1-b)} b^{-k} \varepsilon^{kb} = \varepsilon^k b^{-k}, \quad \lim_{t\to 0^+} h_\varepsilon^k(t) = \lim_{t\to 0^+} (\varepsilon b^{-1})^k = \varepsilon^k b^{-k}$$

Then: $\lim_{t\to 0^-} h_\varepsilon^k(t) = \lim_{t\to 0^+} h_\varepsilon^k(t)$;

when $t=-\varepsilon$,
$$\lim_{t\to -\varepsilon^-} h_\varepsilon^k(t) = 0, \quad \lim_{t\to -\varepsilon^+} h_\varepsilon^k(t) = 0$$

Then: $\lim_{t\to -\varepsilon^-} h_\varepsilon^k(t) = \lim_{t\to -\varepsilon^+} h_\varepsilon^k(t)$.

To sum up, $h_\varepsilon^k(t)$ is continuous in $t=0$, $t=-\varepsilon$, $h_\varepsilon^k(t)$ is continuous in $R$.

(2) To prove $\forall t \in R, h_\varepsilon^m(t) \geq h^m(t)$, Just prove that $\Delta h(t) := h_\varepsilon^k(t) - h^k(t) \geq 0$, By definition:

$$\Delta h(t) := h_\varepsilon^k(t) - h^k(t) := \begin{cases} \Delta h_1(t) = 0, & t \leq -\varepsilon \\ \Delta h_2(t) = \varepsilon^{k(1-b)} b^{-k} (t+\varepsilon)^{kb}, & -\varepsilon < t \leq 0 \\ \Delta h_3(t) = (t+\varepsilon b^{-1})^k - t^k, & t > 0 \end{cases}$$

(i) If $t \leq -\varepsilon$, $\Delta h(t) = \Delta h_1(t) = 0 \geq 0$ is hold;

(ii) If $-\varepsilon < t \leq 0$, then $\varepsilon > 0$, $b > 0$, 且 $t+\varepsilon > 0$, 于是:
$$\Delta h(t) = \Delta h_2(t) = \varepsilon^{k(1-b)} b^{-k} (t+\varepsilon)^{kb} > 0$$

(iii) If $t > 0$, because $\varepsilon > 0$, $b > 0$, then $\varepsilon b^{-1} > 0$:
$$\Delta h(t) = \Delta h_3(t) = (t + \varepsilon b^{-1})^k - t^k > 0$$

To sum up, $\Delta h(t) := h_\varepsilon^k(t) - h^k(t) \geq 0$, then $\forall t \in R, h_\varepsilon^m(t) \geq h^m(t)$

By the construction of $h_\varepsilon^k(t)$,

(i) If $t \leq -\varepsilon$, $\lim_{\varepsilon \to 0} h_\varepsilon^k(t) = 0$ is hold;

(ii) If $-\varepsilon < t \leq 0$, then,
$$\frac{dh_\varepsilon^k(t)}{dt} = k\varepsilon^{k(1-b)} b^{1-k} (t+\varepsilon)^{kb-1},$$

Because $\varepsilon > 0$, $b > 0$, and $t+\varepsilon > 0$, then:
$$\frac{dh_\varepsilon^k(t)}{dt} = k\varepsilon^{k(1-b)} b^{1-k} (t+\varepsilon)^{kb-1} > 0$$

Namely smooth function $h_\varepsilon^k(t)$ is constant monotonous increase in $-\varepsilon < t \leq 0$, then
$$0 = h(-\varepsilon) < h_\varepsilon^k(t) \leq h(0) = \varepsilon^k b^{-k}$$

then
$$0 \leq \lim_{\varepsilon \to 0} h_\varepsilon^k(t) \leq \lim_{\varepsilon \to 0} \varepsilon^k b^{-k} = 0$$

$\lim_{\varepsilon \to 0} h_\varepsilon^k(t) = 0$ is hold.

(iii) If $t > 0$:
$$\lim_{\varepsilon \to 0} h_\varepsilon^k(t) = \lim_{\varepsilon \to 0} (t+\varepsilon b^{-1})^k = t^k = h^k(t)$$

In summary, $\forall t \in R, \lim_{\varepsilon \to 0} h_\varepsilon^m(t) = h^m(t)$

(3) By the construction of $h_\varepsilon^k(t)$, notes:

$$h_\varepsilon^k(t) = \begin{cases} h_1(\varepsilon,t) = 0 & t < -\varepsilon \\ h_2(\varepsilon,t) = \varepsilon^{k(1-b)}b^{-k}(t+\varepsilon)^{kb} & -\varepsilon \leq t < 0 \\ h_3(\varepsilon,t) = (t+\varepsilon b^{-1})^k & t \geq 0 \end{cases}$$

First of all:

$$h_1(t,\varepsilon) = 0 < h_2(t,\varepsilon) \leq \varepsilon^k b^{-k} \leq h_3(t,\varepsilon)$$

Followed by:

$$\frac{dh_2(t,\varepsilon)}{d\varepsilon} = k\varepsilon^{k(1-b)}b^{1-k}(t+\varepsilon)^{kb-1} + k(1-b)\varepsilon^{k(1-b)-1}b^{-k}(t+\varepsilon)^{kb}$$

$$= k\varepsilon^{k(1-b)-1}b^{-k}(t+\varepsilon)^{kb-1}(\varepsilon b + (1-b)(t+\varepsilon))$$

$$= k\varepsilon^{k(1-b)-1}b^{-k}(t+\varepsilon)^{kb-1}(t+\varepsilon - bt)$$

Due to $\varepsilon > 0$, $b > 0$, $-\varepsilon < t \leq 0$, then:

$$\frac{dh_2(t,\varepsilon)}{d\varepsilon} > 0$$

and

$$\frac{dh_3(t,\varepsilon)}{d\varepsilon} = b^{-1}(t+\varepsilon b^{-1})^{k-1} > 0$$

In summary, smooth function $h_\varepsilon^k(t)$ is Constant monotonous increase in $\varepsilon \geq 0$ which means $0 < \varepsilon_1 < \varepsilon_2$, then $\forall t \in R, h_{\varepsilon_1}^k(t) \leq h_{\varepsilon_2}^k(t)$.

**Theorem 2.2.2** if $\frac{1}{b} < k < +\infty$, $\varepsilon > 0$, $h_\varepsilon^k(t)$ is continuous differentiable, and

$$\frac{dh_\varepsilon^k(t)}{dt} = \begin{cases} 0 & t \leq -\varepsilon \\ k\varepsilon^{k(1-b)}b^{1-k}(t+\varepsilon)^{kb-1} & -\varepsilon < t \leq 0 \\ k(t+\varepsilon b^{-1})^{k-1} & t > 0 \end{cases}$$

Consider the low-order penalty function

$$\hat{P}_\varepsilon^k(x,\pi) = \sum_{i=1}^m f_i(x) + \pi\sum_{i=1}^l h_\varepsilon^k(g_i(x)) + \pi\sum_{i=1}^m h_\varepsilon^k(f_i(x) - f_i(\hat{x}))$$

**Theorem 2.2.3** For any $\frac{1}{b} < k < 1$, $\varepsilon > 0$, $b > 0$, and $f_i(x), i \in M$, $g_i(x), i \in I$ is continuous, the following properties hold:

(1) $\hat{P}_\varepsilon^k(x,\pi)$ is continuous in $R^n$;

(2) $\forall x \in R^n, \hat{P}_\varepsilon^k(x,\pi) \geq \hat{P}^k(x,\pi)$;

(3) $\forall x \in R^n, \lim_{\varepsilon \to 0} \hat{P}_\varepsilon^k(x,\pi) = \hat{P}^k(x,\pi)$;

(4) if $0 < \varepsilon_1 \leq \varepsilon_2$, then $\forall x \in R^n, \hat{P}_{\varepsilon_1}^k(x,\pi) \leq \hat{P}_{\varepsilon_2}^k(x,\pi)$, $\forall x \in R^n, \hat{P}_{\varepsilon_1}^k(x,\pi) \leq \hat{P}_{\varepsilon_2}^k(x,\pi)$.

**Theorem 2.2.4** If $\frac{1}{b} < k < 1$, $b > 0$, $\varepsilon > 0$, and $f_i(x), i \in M$, $g_i(x), i \in I$ is continuous differentiable, then $\widehat{P}_\varepsilon^k(x,\pi)$ is continuous differentiable in $R^n$.

$$\nabla_x \widehat{P}_\varepsilon^k(x,\pi) = \sum_{i=1}^{m} \nabla_x f_i(x) + \pi \sum_{i=1}^{l} \nabla_x h_\varepsilon^k(g_i(x)) + \pi \sum_{i=1}^{m} \nabla_x h_\varepsilon^k(f_i(x) - f_i(\hat{x}))$$

# 3 A MOSQP algorithm based on the low-order smooth penalty function

## 3.1 Quadratic programming sub-problem construction

For the initial value diffusion phase, we construct the following problems for each objective function of the multi-objective optimization problem:

$$(\text{P}_i) \quad \min \ f_i(x), i \in M = \{1,\cdots,m\}$$
$$\text{s.t.} \quad g_i(x) \leq 0, i \in I = \{1,2,\cdots a\}$$

The value function we used at this stage is as: $f_i(x)$.

The search direction is obtained by solving the following quadratic planning subproblems:

$$\mathbf{QP}_i(x^k, B^k) \quad \min \frac{1}{2}(d^j)^T B^j d^j + \nabla f_i(x^j)^T d^j, i \in M = \{1,\cdots,m\}$$
$$\text{s.t.} \ \nabla g_i(x^j)^T d + g_i(x^j) \leq 0 \quad i \in I$$

For the multi-objective optimization problem Pareto optimal solution stage, we consider all the objective functions and take each iterative initial value point as the reference point $\hat{x}$, Constructing the problem(P2):

$$(\text{P2}) \quad \min \ \sum_{i=1}^{m} f_i(x)$$
$$\text{s.t.} \quad f_i(x) \leq f_i(\hat{x}), i \in M = \{1,\cdots,m\}$$
$$g_i(x) \leq 0, i \in I = \{1,2,\cdots a\}$$

The corresponding low-order penalty function is:

$$\widehat{P}_\varepsilon^k(x,\pi) = \sum_{i=1}^{m} f_i(x) + \pi \sum_{i=1}^{l} h_\varepsilon^k(g_i(x)) + \pi \sum_{i=1}^{m} h_\varepsilon^k(f_i(x) - f_i(\hat{x}))$$

The search direction is obtained by solving the following quadratic planning subproblems:

$$\mathbf{QP2}(x^k, B^k) \quad \min \frac{1}{2}(d^k)^T B^k d^k + \sum_{i=1}^{m} \nabla f_i(x^k)^T d^k$$

$$s.t. \nabla f_i(x^k)^T d^k + f_i(x^k) - f_i(\hat{x}) \leq 0, i \in M = \{1, \cdots, m\}$$

$$\nabla g_i(x^k)^T d + g_i(x^k) \leq 0, i \in I = \{1, \cdots, l\}$$

## 3.2 A new MOSQP algorithm

**Algorithm 2 (MOSQP algorithm):**

**Iteration point set initialization:** Initialize the iteration point set $X_0 = \{x_1, x_2, \cdots x_N\}, x_i \in R^n$, They are all feasible points.

**Iterative point set feasible spread:** Kth spread for $X_0 = \{x_1, x_2, \cdots x_N\}, x_i \in R^n$.

Step1 Let $j = 0$, $B = I_n$, $0 < \sigma < 1$, $i = 1$, $M > 1$, $0 < A < 1$, $X = X_0$, $T = \emptyset$, $c > 0$;

Step2 Let $x^j \in X_j$, Solving the sub-problem $\mathbf{QP}_i(x^j, I)$, and let $d^j = d(x^j, I_n)$, if $d^j = 0$, $T = T \cup \{x^j\}$, Otherwise turn Step3;

Step3 Let $\alpha^j = 1$, to Step4;

Step4 If $f_i(x^j + \alpha^j d^j) \leq f_i(x^j) + \sigma \alpha^j [\nabla_x f(x^j)]^T d^j$, and $g_i(x^j + \alpha^j d^j) \leq 0, i \in I$, to Step6 Otherwise to Step5;

Step5 Let $\alpha^j = A\alpha^j$, to Step6;

Step6 Let $T = T \cup \{x^j + \alpha^j d^j\}$, $i = i + 1$, if $i = m + 1$ to Step7, Otherwise, to Step2;

Step7 Let $j = j + 1$, $X_j = T$, $i = 0$, $X = X \cup X_j$, if $j = K$, to step8, Otherwise, let $T = \emptyset$, to Step2;

Step8 Output X for the point where the crowding function value $C(x_i)$ is greater than $c$.

**Note 3.2.1** All the points obtained by diffusion steps are in the feasible domain, ensuring that the iteration point enters the Pareto optimal solution stage. When the initial iteration point is taken as the reference point, the feasible domain of the problem (P2) is not empty. Secondly, Step8 is designed to remove repeated points during the diffusion process and close points, and c takes different values according to the order of magnitude corresponding to the specific question.

**Pareto Optimal solution of the multi-objective optimization problem:**

Step1 令 $S = \emptyset$, $\frac{1}{b} < k < 1$, $b > 0$, $\varepsilon > 0$, $j = 0$, $B = I_n$, $0 < \sigma < 1$, $j = 0$, $i = 1$,

$M > 1$, $0 < A < 1$;

Step2 If the subproblem $\mathbf{QP2}(x^j, B^j)$ is feasible, Solving the sub-problem $\mathbf{QP2}(x^j, B^j)$, let $d^j = d(x^j, B^j)$, $d^j = 0$, $S = S \cup x^j$, Otherwise to Step3; If the problem $\mathbf{QP2}(x^j, B^j)$ is not feasible, let $d^j = -\nabla_x P_{\varepsilon^j}^k(x^j, \pi^j)$ to Step8;

Step3 let $\varepsilon^j = A^{k^*} \varepsilon^j$, $st, \varepsilon^j \leq \min_{i \in I_- \neq \emptyset}(-g_i(x^j), f_i(\hat{x}) - f_i(x^j))$,

if $(\nabla_x \hat{P}_{\varepsilon^j}^k(x^j, \pi^j))^T d^j \leq -\frac{1}{2}(d^j)^T B^j d^j$, to Step5, Otherwise to Step4;

Step4 $\pi^j = M\pi^j$, to Step3;

Step5 $\alpha^j = 1$, to Step6;

Step6 if $\hat{P}_{\varepsilon^j}^k(x^j + \alpha^j d^j, \pi^j) \leq \hat{P}_{\varepsilon^j}^k(x^j, \pi) + \sigma \alpha^j [\nabla_x \hat{P}_{\varepsilon^j}^k(x^j, \pi^j)]^T d^j$, to Step13; Otherwise to Step7;

Step7 $\alpha^j = A\alpha^j$, to Step6;

Step8 if $\sum_{g_i(x^j) > 0} d_2(x^j) \nabla_x g_i(x^j) + \sum_{f_i(x^j) > f_i(\hat{x}^j)} d_3(x^j) \nabla_x f_i(x^j) \neq 0$, to step10, Otherwise to Step9;

Step9 let $\varepsilon^j = A\varepsilon^j$, to Step8;

Step10 $(\nabla_x \hat{P}_{\varepsilon^j}^k(x^j, \pi^j))^T d^j \leq -\max\{g_i(x^j), i \in I_+; f_i(x^j) - f_i(\hat{x}), i \in M_+;\} \leq 0$ to Step5, Otherwise to Step12;

Step11 let $\pi^j = M\pi^j$, to Step11;

Step12 If the problem is feasible, Correction of $B^{j+1}$ by using the BFGS formula, Otherwise let $B^{j+1} = B^j$, to Step14;

Step13 let $x^{j+1} = x^j + \alpha^j d^j$, $\varepsilon^{j+1} = \varepsilon^j$, $\pi^{j+1} = \pi^j$, $j = j+1$, $S = S \cup \{x^j + \alpha^j d^j\}$, to Step14;

Step14 let $j = 0$, $i = i+1$, if $i = M+1$, to Step15;

Step15 $S = \text{Non-dominated}(S)$, Output S for the point where the crowding function value $C(x_i)$ is greater than $c$.

**Note3.2.2** The $d_2(x^j)$ and $d_3(x^j)$ in the algorithm are respectively:

$$d_2(x^j) = k(g_i(x^j) + \varepsilon b^{-1})^{k-1}, \quad d_3(x^j) = k(f_i(x^j) - f_i(\hat{x}) + \varepsilon b^{-1})^{k-1}$$

Secondly, the set of the 50 Pareto optimal solutions with the largest crowding degree function value in S for the original multi-objective optimization problem is designed to remove the Pareto front repeat solutions obtained by the MOSQP algorithm.

**Theorem3.2.1** If $x^*$ is the KKT point of the problem (P2), it is the Pareto critical point of the original multi-objective optimization problem.

Proof: If $x^*$ is the KKT point of the problem (P2), then:

$$\sum_{i=1}^{m}\nabla f_i(x^*)+\sum_{i=1}^{m}\mu_i\nabla f_i(x^*)+\sum_{i=1}^{l}\lambda_i^*\nabla g_i(x^*)=0$$

$$g_i(x^*)\leq 0, \lambda_i^*\geq 0, \forall i\in I$$

$$\lambda_i^* g_i(x^*)=0$$

$$f_i(x^*)-f_i(\hat{x})\leq 0, \mu_i\geq 0, i\in M=\{1,\cdots,m\}$$

$$\mu_i(f_i(x^*)-f_i(\hat{x}))=0$$

So $x^*$ satisfied:

$$\sum_{i=1}^{m}(1+\mu_i)\nabla f_i(x^*)+\sum_{i=1}^{l}\lambda_i^*\nabla g_i(x^*)=0$$

$$g_i(x^*)\leq 0, \lambda_i^*\geq 0, \forall i\in I$$

$$\lambda_i^* g_i(x^*)=0$$

So $x^*$ is the KKT point of

$$\min \sum_{i=1}^{m}(1+\mu_i)f_i(x)$$

$$s.t. g_i(x)\leq 0, i\in I$$

And because $1+\mu_i>0$, $x^*$ is the Pareto critical point of the original multi-objective constraint optimization problem (CMOP).

Theorem 3.2.1 also shows that the proposed MOSQP algorithm will converge to the Pareto approximation front of the original multi-objective constraint optimization problem (CMOP).

# 4 Numerical experiments

## 4.1 Algorithm performance performance measure

First, the measure of the algorithm performance is defined as follows:

**Purity metric:** Purity metric is used to compare the number of non-dominated solutions obtained by different algorithms.

$$\bar{t}_{p,s}=\frac{|F_p|}{|F_{p,s}\cap F_p|}$$

From the definition, the smaller $\bar{t}_{p,s}$, the better the quality of the solution. Obviously, $\bar{t}_{p,s}=\infty$ means that the algorithm cannot generate any non-dominant points in the reference Pareto front of the corresponding problem.

**Spread metric:** To analyze whether the points generated by the algorithm are well distributed in the approximate Pareto front of a given problem, two divergence metrics are used ($\Gamma$ and $\Delta_{p,s}$).

$$\Gamma_{p,s} = \max_{j \in M} \max_{i \in \{0,1,\ldots,N\}} \delta_{i,j}$$

$$\Delta_{p,s} = \max_{j \in M} \left( \frac{\delta_{0,j} + \delta_{N,j} + \sum_{i=1}^{N-1} |\delta_{i,j} - \bar{\delta}_j|}{\delta_{0,j} + \delta_{N,j} + (N-1)\bar{\delta}_j} \right)$$

By definition, the smaller the value of the two distributions, the more uniform the Pareto algorithm solves the corresponding problem.

## 4.2 Numerical experiments

The MOSQP algorithm, MOSQP algorithm proposed by Fliege and Vaz in literature [1], and NSGAII algorithm are used for numerical comparison experiments in Python to solve the problem Problem1 ZDT1 ZDT2 MOP3[5-9]. respectively, and the purity index, divergence index and function calculation frequency index are calculated through the algorithm performance measurement formula in Section 4.1. Analysis the quantity and quality of the Pareto, and then compare the performance curve and analyze the advantages and disadvantages of the algorithm. To distinguish between them, we note that the MOSQP algorithm proposed by Fliege and Vaz is MOSQP (F). The numerical experimental parameters of the MOSQP algorithm are configured as follows:

$$k = 0.5, \quad b = 4, \quad \sigma = 0.2, \quad M = 2, \quad A = 0.5, \quad K = 5$$

Where the numerical experimental parameters of the MOSQP (F) algorithm are configured as:

$$\sigma = 0.2, \quad K = 5$$

K is the times of Spread, For the MOSQP algorithm and the MOSQP (F) algorithm for the same initial iteration point set diffusion 5 times. The gradient calculations involved in the algorithm were all solved using the jacobian function in the automated differential library autograd. NSGAII The number of populations is set to the maximum of the number of Pareto fronts acquired by the MOSQP algorithm and the MOSQP (F) algorithm, and the number of iterations is set to 300.

First, we give the following dual-objective constraint optimization problem

(Problem1):
$$\min\ F(x)=[f_1(x)=(1-x_1)^2+100(x_2-x_1^2)^2,$$
$$f_2(x)=(x_1+x_2-1)^2+10(x_1-x_2)^2]^T$$
$$s.t.\quad x_1^2+x_2^2\le 0.5$$

Using Monte Carlo simulations:

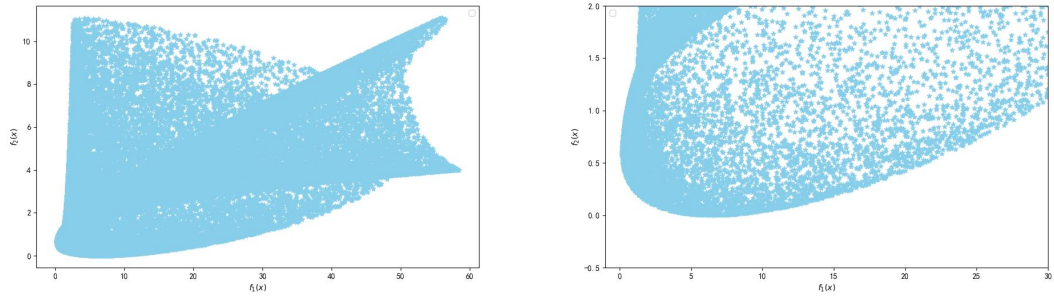

Figure 4,2.1 Functional images generated by the Monte Carlo method

The left side simulates the function image within the entire feasible domain, while the right side simulates the intercepted local features near the Pareto front.
MOSQP, MOSQP (F) and NSGAII algorithms are used to solve the above dual objective constraint optimization problem, and show the solved Pareto frontier:

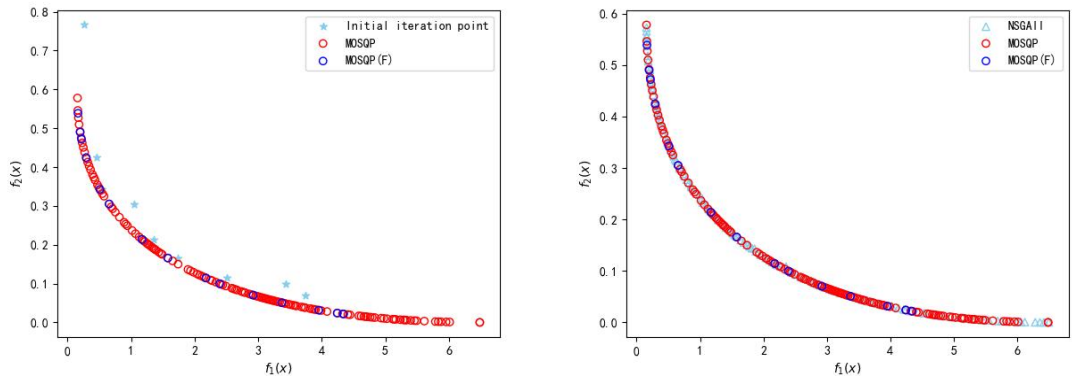

(a) Algorithm iteration process   (b) MOSQP (F) and NSGAII algorithm frontier comparison

Figure 4.2.2. MOSQP algorithm, NSGAII obtained problem (Problem1)

From figure4.2..2, we can see that the MOSQP algorithm proposed in this paper can effectively solve the nonlinear constraint multi-objective optimization problem, and improve the MOSQP algorithm of diffusion step Pareto front in quantity and quality are better than MOSQP (F), we calculated the MOSQP, MOSQP (F) and NSGAII algorithm (300) approximate Pareto front corresponding to the purity metric $\overline{t}_{p,s}$, Spread metric $\Gamma_{p,s}$ $\Delta_{p,s}$.

Table 4.2.1 Solving the corresponding values of each performance index of Problem1

|  | $\overline{t}_{p,s}$ | $\Gamma_{p,s}$ | $\Delta_{p,s}$ |
|---|---|---|---|
| MOSQP | 1.7808 | 0.4737 | 0.8965 |
| MOSQP(F) | 8.6667 | 0.5999 | 0.5438 |
| NSGAII | 3.0952 | 0.2111 | 0.5981 |

From the above indicators, it can be analyzed that based on the same initial iteration point set, the MOSQP algorithm proposed in this paper has more number of non-dominant solutions for the approximate Pareto front obtained by Problem1, with higher quality, which is better than MOSQP (F) and NSGAII algorithm, and the distribution is relatively good. Compared to MOSQP (F), the number of gradient calculation increases, but avoid the calculation of Hessian matrix, the PE index will be smaller.

Then we analyze the classical two-objective boundary constraint optimization problem ZDT 1 (n=30), with 30 decision variables in the problem, which increases the difficulty of the problem, which is defined as follows:

$$(\text{ZDT1}) \quad \min \ F(x) = [f_1(x) = x_1, f_2(x) = 1 - \sqrt{\frac{f_1(x)}{g(x)}}]^T$$

$$g(x) = 1 + \frac{9}{n-1} \sum_{i=2}^{n} x_i$$

$$s.t. \quad 0 \leq x_i \leq 1, i = 1, \cdots, n$$

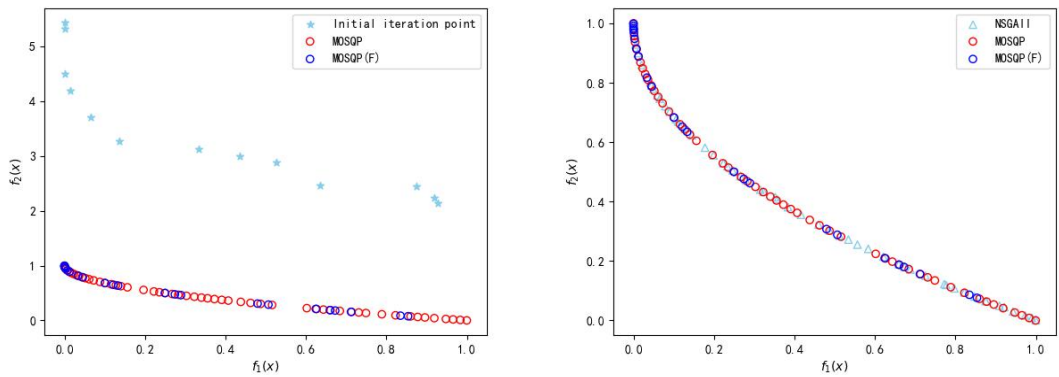

(a) Algorithm iteration process     (b) MOSQP (F) and NSGAII algorithm frontier comparison

Figure 4.2.3 MOSQP algorithm, NSGAII obtained problem (ZDT 1)

From Figure 4.2.3, we can see that the MOSQP algorithm proposed here can effectively solve the high-dimensional boundary constraint multi-objective optimization problem ZDT 1, and improve the MOSQP algorithm of diffusion step Pareto front in quantity and quality are better than MOSQP (F), we calculated the MOSQP, MOSQP (F) and NSGAII algorithm (300) approximate Pareto front corresponding to the purity index $\bar{t}_{p,s}$, Spread index $\Gamma_{p,s}$ $\Delta_{p,s}$.

Table 6.3.2 Solve the corresponding values of each performance index of ZDT 1

|  | $\bar{t}_{p,s}$ | $\Gamma_{p,s}$ | $\Delta_{p,s}$ |
| --- | --- | --- | --- |
| MOSQP | 2.3659 | 0.0996 | 0.5656 |
| MOSQP(F) | 4.2174 | 0.1904 | 0.9252 |
| NSGAII | 2.9394 | 0.66 | 0.4462 |

From the above indicators, it can be analyzed that based on the same initial iteration point set, the MOSQP algorithm proposed in this paper has more number of non-dominated solutions for the approximate Pareto edge obtained by Problem1, with higher quality, which is better than MOSQP (F) and NSGAII algorithms, and the divergence index is also less than MOSQP (F) algorithm. For the problem of high dimension of decision variables, the calculation of Hessian matrix is avoided, and the calculation efficiency can be greatly improved.

ZDT 2 for two-objective boundary constraint optimization of Pareto front (n=30):

$$(\text{ZDT2}) \quad \min \ F(x) = [f_1(x) = x_1, f_2(x) = 1 - (\frac{f_1(x)}{g(x)})^2]^T$$

$$g(x) = 1 + \frac{9}{n-1} \sum_{i=2}^{n} x_i$$

$$s.t. \quad 0 \leq x_i \leq 1, i = 1, \cdots, n$$

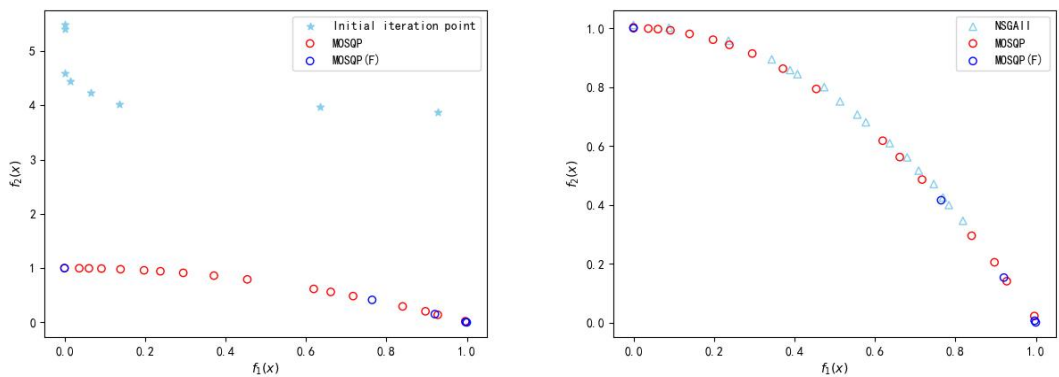

(a) Algorithm iteration process     (b) MOSQP (F) and NSGAII algorithm frontier comparison

Figure 4.2.4 MOSQP algorithm, NSGAII obtained problem (ZDT 2)

From Figure Figure 4.2.4, we can see that the MOSQP algorithm proposed here can effectively solve the frontier non-convex constrained multi-objective optimization problem ZDT 2, and improve the MOSQP algorithm of diffusion step Pareto front in quantity and quality are better than MOSQP (F), we calculated the MOSQP, MOSQP (F) and NSGAII algorithm (300) approximate Pareto front corresponding to the purity metric $\overline{t}_{p,s}$, Spread metric $\Gamma_{p,s}$ $\Delta_{p,s}$.

Table4.2.3 Solve the corresponding values of each performance index of ZDT 2

|          | $\overline{t}_{p,s}$ | $\Gamma_{p,s}$ | $\Delta_{p,s}$ | FE |
|----------|--------|--------|--------|--------|
| MOSQP    | 2.0588 | 0.1906 | 0.6908 | 1148.8 |
| MOSQP(F) | 7      | 0.7643 | 0.7588 | 2219.6 |
| NSGAII   | 2.6923 | 0.1476 | 0.5619 | --     |

From the above indicators, it can be analyzed that based on the same initial iteration point set, the MOSQP algorithm has a higher number of non-dominant solutions for the approximate Pareto edge obtained by Problem1, which is better than MOSQP (F) and NSGAII algorithms, and the divergence indexes are less than MOSQP (F) algorithm. However, the divergence index will deviate relative to the NSGAII algorithm. However, the divergence index will be higher relative to the NSGAII algorithm.

For the two-objective boundary constraint optimization problem of the Pareto front MOP 3:

$$(\text{MOP3}) \quad \max \ F(x) = [f_1(x) = -1 - (A_1 - B_1)^2 - (A_2 - B_2)^2,$$
$$f_2(x) = -(x_1 + 3)^2 - (x_2 + 1)^2]^T$$
$$A_1 = \sin 1 - 2\cos 1 + \sin 2 - 1.5\cos 2$$
$$A_2 = 1.5\sin 1 - \cos 1 + 2\sin 2 - 0.5\cos 2$$
$$B_1 = \sin x_1 - 2\cos x_1 + \sin x_2 - 1.5\cos x_2$$
$$B_2 = 1.5\sin x_1 - \cos x_1 + 2\sin x_2 - 0.5\cos x_2$$
$$\text{s.t.} \quad -\pi \leq x_i \leq \pi, i = 1, \cdots, n$$

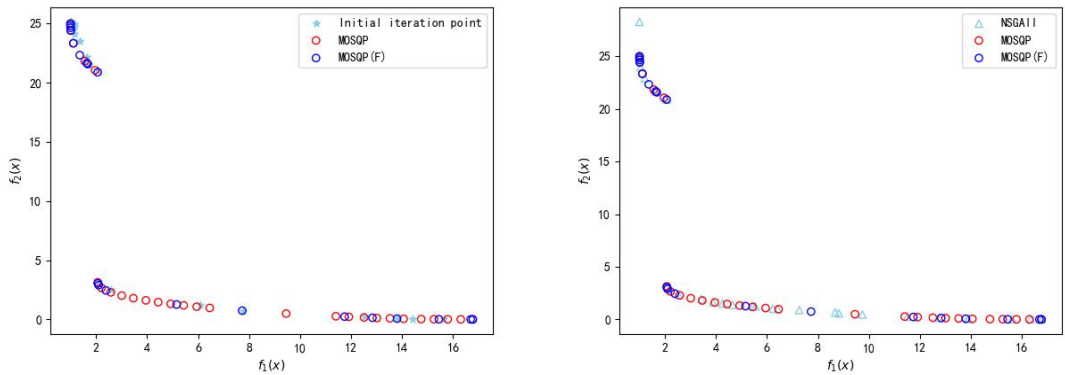

(a) Algorithm iteration process    (b) MOSQP (F) and NSGAII algorithm frontier comparison

Figure 4.2.5. MOSQP algorithm, NSGAII obtained problem (MOP3)

From Figure 4.2.5, we can see that the MOSQP algorithm proposed in this paper can effectively solve the noncontinuous constrained multi-objective optimization problem of Pareto front, MOP 3, and in this paper propose that the Pareto front obtained by MOSQP algorithm is better than MOSQP (F) algorithm in both quantity and quality.same. we calculated the MOSQP, MOSQP (F) and NSGAII algorithm (300) approximate Pareto front corresponding to the purity metric $\bar{t}_{p,s}$, Spread metric $\Gamma_{p,s}$ $\Delta_{p,s}$.

Table 4.2.4 Solving the corresponding values of each performance index of MOP 3

|  | $\bar{t}_{p,s}$ | $\Gamma_{p,s}$ | $\Delta_{p,s}$ |
| --- | --- | --- | --- |
| MOSQP | 2.1228 | 17.7978 | 1.4662 |
| MOSQP(F) | 6.05 | 17.8549 | 1.2833 |
| NSGAII | 2.75 | 17.7383 | 1.4027 |

It can be analyzed from Table 4.2.4 that, based on the same initial iteration point set, the proposed MOSQP algorithm has a larger number of non-dominant solutions for the approximate Pareto front obtained by Problem1, with higher quality, which is better than MOSQP (F) and NSGAII algorithm, and the divergence index is not much different from MOSQP (F) and NSGAII algorithm.

By solving the above different types of problems, we find that the MOSQP algorithm proposed can be effectively solved for nonlinear function constraint optimization, high-dimensional boundary constraint optimization, non-convex in Pareto front, and non-continuous multi-objective optimization in Pareto front. And the improved MOSQP is better than the MOSQP (F) algorithm and NSGAII in the number of non-dominant solutions in the Pareto front. The divergence index is inferior to NSGAII, but the MOSQP (F) algorithm has different problems. Moreover, in the multi-objective optimization problem with high dimension of decision variables, it can improve the computational efficiency.

# 5 Conclusion

This paper proposes a new MOSQP algorithm based on the low-order smooth penalty function as the value function of line search. The algorithm constructs a single objective optimization problem based on each objective function, and solves the secondary planning sub-problem to spread in the feasible domain, selects the non-dominant point set after diffusion, then constructs the new secondary planning problem by Armijo step rule, linear search step, and corrected by Powell in 1978. If the secondary planning problem fails, the negative gradient of the value function is used as the search direction. We also show that the algorithm converges to the approximate Pareto frontier for the multiobjective constraint optimization problem under several basic assumptions. By solving the above different types of problems, we find that the MOSQP algorithm proposed can be effectively solved for nonlinear function constraint optimization, high-dimensional boundary constraint optimization, non-convex in Pareto front, and non-continuous multi-objective optimization in Pareto front. And the improved MOSQP is better than the MOSQP (F) algorithm and NSGAII in the number of non-dominant solutions in the Pareto front. The divergence index is inferior to NSGAII, but the MOSQP (F) algorithm has different problems.